\documentclass[12pt,letterpaper]{article}
\usepackage{amsfonts}
\usepackage{amssymb}
\usepackage{latexsym}  
\usepackage{pifont}  
\usepackage{bbm}  %

\let\SavedRightarrow=\Rightarrow
\usepackage{marvosym}
\let\Rightarrow=\SavedRightarrow

\usepackage{amsmath}   

\newcommand\QQQ{{\mathbb Q}}
\newcommand\PPP{{\mathbb P}}
\newcommand\FFF{{\mathbb F}}
\newcommand\RRR{{\mathbb R}}

\newcommand\cccc{{\mathfrak c}}

\newcommand\V{{\mathbf V}}  

\newcommand\one{\mathbbm{1}} 

\newcommand\cl{\mathrm{cl}}   
\newcommand\diam{\mathrm{diam}}   

\newcommand{\forces}{\Vdash}

\newcommand\down{\mathord {\downarrow}}  

\newcommand\hgt{\mathrm{ht}}   
\newcommand\MA{\mathrm{MA}}  
\newcommand\CH{\mathrm{CH}}  

\newcommand\iv{^{-1}} 

\newcommand\onto{\twoheadrightarrow}

\newcommand\eop{{\Large \Coffeecup}}  

\newenvironment{itemizz}{\begin{itemize}\setlength{\itemsep}{-1mm}} %
{\end{itemize}}                              

\newenvironment{itemizn}[1] 
{\begin{itemize} \setlength{\itemsep}{-1mm} %
\renewcommand\labelitemi{\ding{#1}}} %
{\end{itemize}}

\newtheorem{theorem}{Theorem}[section]
\newtheorem{definition}[theorem]{Definition}
\newtheorem{lemma}[theorem]{Lemma}

\newenvironment{proof}{{\bf Proof.}}{\eop\medskip}
\newenvironment{proofof}[1]{\medskip \textbf{Proof of #1.}}{\eop\medskip}

\begin{document}

\title{Locally Connected  HL Compacta}

\author{Kenneth Kunen\footnote{University of Wisconsin,  Madison, WI  53706, U.S.A.,
\ \ kunen@math.wisc.edu}
\thanks{Author partially supported by NSF Grant DMS-0456653.}  }

\maketitle

\begin{abstract}
It is consistent with $\MA + \neg \CH$ that there is a locally
connected hereditarily Lindel\"of compact space which is not metrizable.
\end{abstract}

\section{Introduction} 
\label{sec-intro}
All spaces discussed in this paper are assumed to be Hausdorff.
A question attributed in 1982 by Nyikos \cite{MER} to M. E. Rudin
asks whether $\MA + \neg\CH$ implies that every locally
connected hereditarily Lindel\"of (HL) compact space is metrizable
(equivalently, second countable); see Gruenhage \cite{Gr} for further
discussion.
Filippov \cite{Fil} had constructed such a space
in 1969 under CH, and his space is also hereditarily separable (HS).
Since Filippov used a Luzin set in his construction,
and $\MA + \neg \CH$ implies that there are no Luzin sets,
it might have been hoped that $\MA + \neg \CH$ refutes the existence
of such a space, but that turns out to be false; we shall show in Section
\ref{sec-oca-fail}:

\begin{theorem}
\label{thm-main}
It is consistent with $\MA + 2^{\aleph_0} = \aleph_2$ that
there is a non-metrizable locally connected compactum
which is both HS and HL.
\end{theorem}

Our proof shows in ZFC that the Filippov construction succeeds
provided that there is a \emph{weakly Luzin set}; details
are in Section \ref{sec-wk-luzin}.
Weakly Luzin sets are related to entangled sets, and our proof of Theorem \ref{thm-main} 
shows that weakly Luzin sets are consistent
with $\MA + 2^{\aleph_0} = \aleph_2$.
We can show that PFA refutes spaces which are ``like'' the Filippov
space (see Section \ref{sec-SOCA}),
but we do not know whether PFA refutes all
non-metrizable locally connected HL compacta.

The Filippov space may be viewed as a connected version of
the double arrow space $D$, which was described in 1929 by
Alexandroff and Urysohn \cite{AU}.  This is a ZFC example
of a non-metrizable compactum which is both HS and HL,
but it is totally disconnected.  The cone over $D$ yields a connected
example, but this is not locally connected.

$D$ is constructed from $[0,1]$ by replacing the points of $(0,1)$
by neighboring pairs of points.  To construct the Filippov space,
start with $[0,1]^2$,
choose a set $E \subseteq (0,1)^2$, and replace the points of $E$ by circles,
obtaining a space $\Phi_E$.
This $\Phi_E$ is compact and locally connected.
$\Phi_E$ is metrizable iff $E$ is countable.
Furthermore, if $E$ is a Luzin set, then, as
Filippov showed, $\Phi_E$ is HL, and a similar proof shows that
$\Phi_E$ is HS as well.

Actually, by Juh\'asz \cite{Ju} and Szentmikl\'ossy \cite{Sze},
HS and HL are equivalent for compacta under $\MA(\aleph_1)$,
but that result is not needed here.  We shall show in ZFC
(Theorem \ref{thm-equivs})
that $\Phi_E$ is HS iff $\Phi_E$ is HL iff $E$ is weakly Luzin.

\section{Weakly Luzin Sets}
\label{sec-wk-luzin}
We begin by describing Filippov's example \cite{Fil}.
We start with $[0,1]^n$ (where $1 \le n < \omega$),
rather than $[0,1]^2$, to show that
the construction does not depend on accidental features
of two-dimensional geometry.
As usual, $S^{n-1} \subset \RRR^n$ denotes the unit sphere,
and $\| x \|$ denotes the length of $x \in \RRR^n$, using the
standard Pythagorean metric.
Given $E \subseteq (0,1)^n$, we shall obtain the space
$\Phi_E$ by replacing all points in $E$ by $(n-1)$--spheres and leaving
the points in $[0,1]^n \backslash E$ alone.

\begin{definition}
$\rho : \RRR^n \backslash \{0\} \onto S^{n-1}$
is the perpendicular retraction: $\rho(x) = x/ \|x\|$.
\end{definition}

So, $\rho(y-x)$ may be viewed as the direction from $x$ to $y$.

\begin{definition}
Fix $E \subseteq (0,1)^n$
and let $E' = [0,1]^n \backslash E$.
The \emph{Filippov space} $\Phi_E$, as a set,
is $(E \times S^{n-1}) \cup E'$.
Define $\pi = \pi_E : \Phi_E \onto [0,1]^n $ so that
$\pi(x,w) = x$ for $(x,w) \in E \times S^{n-1}$, and
$\pi(x) = x$ for $x \in E'$.
For $\varepsilon > 0$,  define, for $x \in E'$:
\[
B(x, \varepsilon) = \{p \in \Phi_E : \|\pi(p) - x\| < \varepsilon\}\ \ ,
\]
and define, for $x \in E $ and $W$ an open subset of $S^{n-1}$:
\[
B(x,W,\varepsilon) = \{x\} \times W
\ \cup \ \{p \in \Phi_E :0 < \|\pi(p) - x\| < \varepsilon 
\ \& \ \rho(\pi(p) -  x) \in W \}\ \ .
\]
Give $\Phi_E$ the topology which has all the sets
$B(x, \varepsilon)$  and $B(x,W,\varepsilon)$ as a base.
\end{definition}

\begin{lemma}
For each $E \subseteq (0,1)^n$:
$\Phi_E$ is compact and first countable.   $\pi_E$ is 
a continuous irreducible map from $\Phi_E$ onto $[0,1]^n$.
$\Phi_E$ is metrizable iff $E$ is countable.
If $n \ge 2$, then $\Phi_E$ is connected and
locally connected, and $\pi_{E}$ is monotone.
\end{lemma}

The proof of this last sentence uses the connectedness of $S^{n-1}$.
When $n = 1$,  $S^0 = \{\pm 1\}$, and
$\Phi_E$ is just the double arrow space obtained
by doubling the points of $E$, so $\Phi_E$ 
is always HS and HL.
When $n > 1$, the argument of Filippov shows that
$\Phi_E$ is HL if $E$ is a Luzin set,
but actually something weaker than Luzin suffices:

\begin{definition} 
For $1 \le n < \omega$:
\begin{itemizn}{42}
\item
If $T \subseteq \RRR^n$, then $T^* = \{x-y : x,y \in T \ \&\ x \ne y\}$
\item
$T \subseteq \RRR^n$ is \emph{skinny} iff $\cl(\rho(T^*)) \ne S^{n-1}$.
\item
$E \subseteq \RRR^n$ is a \emph{weakly Luzin} set iff
$E$ is uncountable and every skinny subset of $E$ is countable.
\end{itemizn}
\end{definition}

Every subset of a skinny set is skinny, and $T$ is skinny
iff $\overline T$ is skinny.  Each skinny set is nowhere dense,
so every Luzin set is weakly Luzin.
When $n = 1$,  $T$ is skinny iff $|T| \le 1$,
every uncountable set is weakly Luzin,
and the proof of the following theorem reduces to the usual proof that the
double arrow space is HS and HL.

When $n > 1$: Under CH, it is easy to construct a weakly Luzin set
which is not Luzin.  PFA implies that there are no weakly Luzin sets.
We shall show in Section \ref{sec-oca-fail} that a weakly Luzin set is
consistent with $\MA + \cccc = \aleph_2$.
Clearly, if there is a weakly Luzin set in $\RRR^n$, then there
is one in $(0,1)^n$.

\begin{theorem}
\label{thm-equivs}
For $n \ge 1$ and uncountable $E \subseteq (0,1)^n$,
the following are equivalent:
\begin{itemizz}
\item[1.] $E$ is weakly Luzin.
\item[2.] $\Phi_E$ is HS.
\item[3.] $\Phi_E$ is HL.
\item[4.] $\Phi_E$ has no uncountable discrete subsets.
\end{itemizz}
\end{theorem}
\begin{proof}
For $(4) \to (1)$:  If $E$ is not weakly Luzin, fix an uncountable
skinny $T \subseteq E$.  Let $W = S^{n-1} \backslash \cl(\rho(T^*))$,
and fix $w \in W$.  Then $\{(x,w) : x \in T\} \subset \Phi_E$ is 
discrete.

Since $(2) \to (4)$ and $(3) \to (4)$ are obvious, it is sufficient
to prove $(1) \to (2)$ and $(1) \to (3)$.  So, assume (1),
and let $\langle p_\alpha : \alpha < \omega_1 \rangle$ be an
$\omega_1$--sequence of distinct points from $\Phi_E$; we show that
it can be neither left separated nor right separated.
To do this, fix an open neighborhood $N_\alpha$ of $p_\alpha$ for each 
$\alpha$; we find $\alpha < \beta < \gamma$ such that 
$p_\beta \in N_\alpha$ and $p_\beta \in N_\gamma$.
This is trivial if $\aleph_1$ of the $\pi(p_\alpha)$ lie in $E'$,
or if $\aleph_1$ of the $\pi(p_\alpha)$ are the same point of $E$.
So, thinning the sequence (discarding some points),
and shrinking the neighborhoods (replacing them by smaller ones),
we may assume that
each $p_\alpha = (x_\alpha, w_\alpha) \in E \times S^{n-1}$ and that
$N_\alpha = B(x_\alpha, W, \varepsilon)$, where the
$x_\alpha$ are distinct points in $E$,
$W$ is open in $S^{n-1}$, and each $w_\alpha \in W$.
Let $T = \{x_\alpha : \alpha < \omega_1\}$.  Thinning further,
we may assume that $\diam(T) < \varepsilon$, so that
$p_\beta \in N_\alpha$ iff $\rho(x_\beta - x_\alpha) \in W$.
Thinning again, we may assume that every point of $T$ is
a condensation point of $T$.
Since $E$ is weakly Luzin, $T$ cannot be skinny, so
$\rho(T^*)$ is dense in $S^{n-1}$, so fix $\xi \ne \eta$
such that $\rho(x_\eta - x_\xi) \in W$.
There are then open $U \ni x_\xi$ and $V \ni x_\eta$ such
that $\rho(z - y) \in W$ for all $y \in U$ and $z \in V$.
Since $|U \cap T| = |V \cap T| = \aleph_1$, we may fix
$\alpha < \beta < \gamma$ with $x_\alpha,x_\gamma \in U$ and $x_\beta \in V$;
then
$\rho(x_\beta - x_\alpha) \in W$ and
$\rho(x_\beta - x_\gamma) \in W$, so
$p_\beta \in N_\alpha$ and $p_\beta \in N_\gamma$.
\end{proof}

Entangled subsets of $\RRR$ were discussed by Avraham and Shelah \cite{AS}
(see also \cite{ARS}).
The weakly Luzin sets and the entangled sets have a common generalization:

\begin{definition} 
For $1 \le n < \omega$ and $1 \le k < \omega$:
\begin{itemizz}
\item[1.]
If $E \subseteq \RRR^n$, then $\widetilde E \subseteq (\RRR^n)^k$
is \emph{derived} from $E$ iff $\widetilde E \subseteq E^k$ and
whenever $\vec x = \langle x_0, \ldots x_{k-1} \rangle \in \widetilde E$
and $\vec y = \langle y_0, \ldots y_{k-1} \rangle \in \widetilde E$:
$x_i \ne y_j$ unless $i = j$ and $\vec x = \vec y$.
\item[2.]
$E$ is $(n,k)$--\emph{entangled} iff 
$E \subseteq \RRR^n$ is uncountable and whenever 
$\widetilde E \subseteq (\RRR^n)^k$ is uncountable and derived from $E$,
and, for $i < k$, $W_i$ is open in $S^{n-1}$ with $W_i \ne \emptyset$:
there exist $\vec x, \vec y \in \widetilde E$ with 
$\vec x \ne \vec y$ and $\rho(x_i - y_i) \in W_i$ for all $i$.
\end{itemizz}
\end{definition}

Then ``weakly Luzin''  is equivalent to ``$(n,1)$--entangled'', and
``$k$--entangled'' is equivalent to ``$(1,k)$--entangled''.
$E \subseteq \RRR$ is $(1,1)$--entangled iff $E$ is uncountable.
If $E$ is $(n,k)$--entangled and 
$\widetilde E $ and the $W_{i}$ are as in (2),
then there are actually uncountable disjoint
$X,Y \subseteq \widetilde E  $ such
that $\forall i\, \rho(x_i - y_i) \in W_i$
whenever $\vec x \in X$ and $\vec y \in Y$.
In (2), when $k=1$, WLOG we may assume that $W_0= -W_0 $.

\section{Preserving Failures of SOCA}
\label{sec-oca-fail}
The \emph{Semi Open Coloring Axiom} (SOCA) is a well-known consequence of the PFA;
see Abraham, Rubin, and Shelah \cite{ARS}.
We shall show that certain classes
of failures of SOCA can be preserved in ccc extensions
satisfying $\MA + 2^{\aleph_0} = \aleph_2$.
This is patterned after the proof (see \cite{ARS, AS})
that an entangled set is consistent
with $\MA + 2^{\aleph_0} = \aleph_2$.

\begin{definition}
For any set $E$:
Let $E^\dag = (E \times E) \setminus \{(x,x) : x \in E\}$.
Fix $W \subseteq E^\dag$ with $W = W\iv$.  Then
$T \subseteq E$ is \emph{$W$--free} iff
$T^\dag \cap W = \emptyset$ and
$T$ is \emph{$W$--connected} iff $T^\dag \subseteq W$.
\end{definition}

\begin{definition}
$(E,W)$ is \emph{good} iff $E$ is an uncountable separable metric space,
$W = W\iv$ is an open subset of $E^\dag$, and no uncountable subset of $E$
is $W$--free.
\end{definition}

Then, the SOCA is the assertion that whenever $(E,W)$ is good,
there is an uncountable $W$--connected set.
An uncountable
$E \subseteq \RRR^n$ is weakly Luzin iff $(E,W)$ is good for all
$W$ of the form $\{(x,y)\in E^\dag : \rho(x-y) \in A\}$,
where $A \subseteq S^{n-1}$ is open and $A = -A \ne \emptyset$.
We shall prove:

\begin{theorem}
\label{thm-oca-fail}
Assume that in the ground model $\V$, 
$\CH + 2^{\aleph_1} = \aleph_2$ holds and $E$ is a separable metric space.
Then there is a ccc extension
$\V[G]$ satisfying $\MA + 2^{\aleph_0} = \aleph_2$ such that
for all $W \in \V$, if
$(E,W)$ is good in $\V$ then $(E,W)$ is good in $\V[G]$.
\end{theorem}

A good $(E,W)$ does not by itself contradict SOCA,
since there may be an uncountable subset of $E$ which
is $W$--connected.  But, if $(E,U)$ and $(E,W)$ are both
good and $U \cap W = \emptyset$, then SOCA is contradicted,
since any $W$--connected set is $U$--free.
Such $E,U,W$ are provided by a weakly Luzin
$E \subseteq \RRR^n$ (for $n \ge 2$).
The following combinatorial lemma will be used in the proof of Theorem \ref{thm-oca-fail}.

\begin{lemma}
\label{lemma-comb-good}
Assume the following:
\begin{itemizz}
\item[1.]
CH holds.
\item[2.]
$m \in \omega$; and $(E, W_i)$ is good for each $i \le m$.
\item[3.]
$\theta$ is a suitably large regular cardinal and
$\langle M_\xi : \xi < \omega_1 \rangle$ is a continuous 
chain of countable elementary submodels of $H(\theta)$,
with $E \in M_0$ and each $M_\xi \in M_{\xi + 1}$.
\item[4.]
For $x \in \bigcup_\xi M_\xi \setminus M_0$:
$\hgt(x)$ is the $\xi$ such that $x \in M_{\xi+1} \backslash M_\xi$.
\item[5.]
$ x_\alpha^i \in E \backslash M_0$ for
$\alpha < \omega_1$ and $i \le m$.
\item[6.]
$\hgt( x_\alpha^i) \ne \hgt( x_\beta^j)$ unless
$\alpha = \beta$ and $i = j$.
\end{itemizz}
Then there are $\alpha \ne \beta$ such that 
$(x^i_\alpha , x^i_\beta) \in W_i$ for all $i$.
\end{lemma}

We remark that (6) expresses the standard trick of using a set of points
spaced by a chain of elementary submodels.
In (5), we say $ x_\alpha^i \in E  \backslash M_0$
so that $\hgt( x_\alpha^i)$ is defined; note that by CH,
$ E \subset \bigcup_\xi M_\xi$.

\begin{proof}
Induct on $m$.  When $m = 0$, this is immediate from
the fact that $(E, W_0)$ is good.
Now, assume the lemma for $m-1$, and we prove it for $m$.
Let $\vec x_\alpha= \langle x_\alpha^0, \ldots, x_\alpha^m\rangle \in E^{m+1}$.
Let $\xi(\alpha,i) = \hgt( x_\alpha^i)$.  Thinning 
the $\omega_1$--sequence and rearranging each $\vec x_\alpha$ if necessary,
we may assume that 
$\xi(\alpha, 0) < \xi(\alpha, 1) < \cdots < \xi(\alpha, m)$ and that
$\alpha < \beta \to \xi(\alpha, m) < \xi(\beta, 0)$.  Let
$F = \cl\{\vec x_\alpha : \alpha < \omega_1\} \subseteq  E^{m+1}  $,
and fix
$\mu < \omega_1$ such that $F \in M_\mu$; there is such a $\mu$ by CH.

For $\alpha \ge \mu$:  Let $K_\alpha =
\{z \in E : \langle x_\alpha^0, \ldots, x_\alpha^{m-1}, z\rangle \in F\}$.
$K_\alpha$ is uncountable because 
$K_\alpha \in M_{\xi(\alpha,m)}$ but $K_\alpha$ contains the
element $x_\alpha^m \notin M_{\xi(\alpha,m)}$.
Since $(E, W_m)$ is good, choose $u_\alpha,v_\alpha \in K_\alpha$
with $(u_\alpha , v_\alpha) \in W_m$,
and then choose disjoint basic open sets $U_m, V_m \subseteq E$ with 
$u_\alpha \in U_m$, $v_\alpha \in V_m$,  and 
$(x,y) \in W_m$ for all $x \in U_m$ and $y \in V_m$.

Of course, $U_m, V_m$ depend on $\alpha$, but there are only
$\aleph_0$ possible choices, so fix an uncountable set
$I \subseteq \{\alpha : \mu \le \alpha < \omega_1\}$ such that
the $U_m, V_m$ are the same for $\alpha \in I$.
By the lemma for $m-1$, fix $\gamma,\delta \in I$ such that
$\gamma \ne \delta$ and
$(x^i_\gamma , x^i_\delta) \in W_i$ for all $i < m$.
Now choose disjoint open neighborhoods
$U_i$ of $x_\gamma^i$ and $V_i$ of $x_\delta^i$ for $i < m$
so that $(x , y) \in W_i$ whenever $x \in U_i$ and $y \in V_i$.
Note that the two open sets $\prod_{i \le m} U_i$ and
$\prod_{i \le m} V_i$ both meet $F$,
since $u_\gamma \in K_\gamma$ and $v_\delta \in K_\delta$,
so
$\langle x_\gamma^0, \ldots, x_\gamma^{m-1}, u_\gamma \rangle \in F
\cap \prod_{i \le m} U_i$ and
$\langle x_\delta^0, \ldots, x_\delta^{m-1}, v_\delta \rangle \in F
\cap \prod_{i \le m} V_i$.
We may then choose $\alpha,\beta$ such that
$\vec x_\alpha \in \prod_{i \le m} U_i$ and
$\vec x_\beta \in \prod_{i \le m} V_i$.  But then
$(x^i_\alpha , x^i_\beta) \in W_i$ for all $i$.
\end{proof}

\begin{lemma}
\label{lemma-kill-ccc}
In the ground model $\V$: Assume CH, let $(E,W)$  be good, and let $\QQQ$
be any forcing poset such that $q \forces_\QQQ$ ``$(E,W)$ is not good''
for some $q \in \QQQ$.

Then, in $\V$: there is a ccc poset $\PPP$
of size $\aleph_1$ such that $\QQQ \times \PPP$ is
not ccc and such that for \emph{all} $U \in \V$:
If $(E,U)$ is good then $\one \forces_\PPP$ ``$(E,U)$ is good''.
\end{lemma}
\begin{proof}
Extending $q$, we may assume that for some $\QQQ$--name 
$\mathring Z$: $q \forces  $ ``$\mathring Z \subseteq E$ is uncountable and
$W$--free''.
Fix $\theta$ and the $M_\xi$ so that (3)(4) of Lemma \ref{lemma-comb-good} hold.

Now, inductively choose $q_\alpha \le q$ and
$x_\alpha^0, x_\alpha^1 \in E \backslash M_0$ for $\alpha < \omega_1$
so that $q_\alpha \forces x_\alpha^0, x_\alpha^1 \in \mathring Z$
and such that $\hgt(x_\alpha^0)<\hgt(x_\alpha^1)<\hgt(x_\beta^0)$
whenever $\alpha < \beta < \omega_1$.
Let 
\[
\PPP = \left\{p \in [\omega_1]^{< \omega}  :
\forall \{\alpha,\beta\} \in [p]^2 \,
\left[(x^0_\alpha , x^0_\beta) \in W \text{ or }
(x^1_\alpha , x^1_\beta) \in W \right]\right\} \ \ .
\]
$\PPP$ is ordered by reverse inclusion, with $\one = \emptyset$.
Each $\{\alpha\} \in \PPP$, and the pairs
$(q_\alpha, \{\alpha\}) \in \QQQ \times \PPP$ are incompatible,
so $\QQQ \times \PPP$ is not ccc.

Now, suppose that we have some good
$(E,U)$ and $p \forces_\PPP$ ``$(E,U)$ is not good'';
we shall derive a contradiction.
Extending $p$, we may assume that for some $\PPP$--name 
$\mathring T$: $p \forces  $ ``$\mathring T \subseteq E$ is uncountable and
$U$--free''.
Then, inductively choose $p_\mu \le p$ and
$t_\mu \in E \backslash M_0$ for $\mu < \omega_1$
so that $p_\mu \forces t_\mu \in \mathring T$ and
such that $\hgt(t_\mu)<\hgt(t_\nu)$
whenever $\mu < \nu < \omega_1$.
Our contradiction will use the observation:
\[
\mu \ne \nu \to
(t_\mu, t_\nu) \notin U \text{ or }
p_\mu \perp p_\nu  \ \ . \tag{$*$}
\]
Thinning the sequence and extending $p$ if necessary, we may assume that
the $p_\mu$ form a $\Delta$ system with root $p$; so
$p_\mu =  p \cup \{\alpha(0,\mu), \ldots, \alpha(c,\mu)\}$, with
$\alpha(0,\mu) < \ldots < \alpha(c,\mu)$.
We also assume that $\max(p) < \alpha(0,0)$ and
$\mu < \nu \to \alpha(c,\mu) < \alpha(0,\nu)$.
Since $p_\mu \in \PPP$, 
\[
i \ne j \to
(x^0_{\alpha(i, \mu)} , x^0_{\alpha(j, \mu)}) \in W \text{ or }
(x^1_{\alpha(i, \mu)} , x^1_{\alpha(j, \mu)}) \in W 
\]
for each $\mu$.  Let $\vec x_\mu = 
(x^0_{\alpha(0, \mu)} , x^1_{\alpha(0, \mu)} \ldots
x^0_{\alpha(c, \mu)} , x^1_{\alpha(c, \mu)}) \in E^{2(c+1)} $.
Since $W$ is open, we may thin again and assume that all $\vec x_\mu$
are sufficiently close to some condensation point
of $\{\vec x_\mu : \mu < \omega_1\}$ so that for all $\mu,\nu$:
\[
i \ne j \to
(x^0_{\alpha(i, \mu)} , x^0_{\alpha(j, \nu)}) \in W \text{ or }
(x^1_{\alpha(i, \mu)} , x^1_{\alpha(j, \nu)}) \in W  \ \ .
\]
Thus, if $p_\mu \perp p_\nu$ then the incompatibility must come
from the same index $i$, so that $(*)$ becomes
\[
\mu \ne \nu \to
(t_\mu, t_\nu) \notin U \text{ or }  \exists i\le c \, \left[
(x^0_{\alpha(i, \mu)} , x^0_{\alpha(i, \nu)}) \notin W \text{ and }
 (x^1_{\alpha(i, \mu)} , x^1_{\alpha(i, \nu)}) \notin W \right]  \ \ .
\]
This comes close to contradicting Lemma \ref{lemma-comb-good}.
With an eye to satisfying hypothesis (6), we thin the sequence again
and assume that $\hgt(t_\mu) \ne \hgt(x^\ell_{\alpha(i, \nu)})$
whenever $\mu \ne \nu$.
It is still possible to have $\hgt(t_\mu) = \hgt(x^\ell_{\alpha(i, \mu)})$,
but for each $\mu$, $\hgt(t_\mu) = \hgt(x^\ell_{\alpha(i, \mu)})$
can hold for at most one pair $(\ell,i)$.
Thinning once more, we can assume WLOG that this $\ell$ is always $1$,
so that $\hgt(t_\mu) \ne \hgt(x^0_{\alpha(i, \nu)})$
for all $\mu < \omega_1$ and all $i \le c$.
But now the $(c+2)$--tuples
$(t_\mu,  x^0_{\alpha(0, \mu)}, \ldots  x^0_{\alpha(c, \mu)})$ (for $\mu < \omega_{1} $)
contradict Lemma \ref{lemma-comb-good}, where $W_0 = U$
and the other $W_i = W$.

We also need to show that $\PPP$ is ccc.  If this fails,
then choose the $p_\mu$  to enumerate an antichain.
Derive a contradiction as before,
but replace $(*)$  by the stronger fact
$\mu \ne \nu \to p_\mu \perp p_\nu$,
and delete all mention of $\mathring T$ and the $t_\mu$.
\end{proof}

We remark that a simplification of the above proof yields
the standard proof that an instance of SOCA
can be forced by a ccc poset.  Forget about $\QQQ$ and
just assume that $(E,W)$ is good.  Choose the
$x_\alpha \in E \backslash M_0$ for $\alpha < \omega_1$
so that $\hgt(x_\alpha)<\hgt(x_\beta)$ whenever $\alpha < \beta < \omega_1$.
$\PPP$ is now $\{p \in [\omega_1]^{< \omega}  :
\forall \{\alpha,\beta\} \in [p]^2 \,
[(x_\alpha , x_\beta) \in W \text ]\}$.
Then some $p \in \PPP$ forces an uncountable $W$--connected set.

\begin{proofof}{Theorem \ref{thm-oca-fail}}
In the ground model $\V$, 
we build a normal chain of ccc posets,
$\langle \FFF_\alpha : \alpha \le \omega_2 \rangle$,
where $\alpha < \beta \to \FFF_\alpha \subseteq_c \FFF_\beta$ and we take 
unions at limits.  So, our model will be the generic extension $\V[G]$
given by $\FFF_{\omega_2}$.  
$|\FFF_\alpha| \le \aleph_1$ for all $\alpha < \omega_2$, while
$|\FFF_{\omega_2}| = \aleph_2$.  
Given $\FFF_\alpha$, we choose
$\mathring{\PPP}_\alpha$, which is an $\FFF_\alpha$--name 
forced by $\one$ to be a ccc\@ poset of size $\aleph_1$;
then $\FFF_{\alpha+1} = \FFF_\alpha *  \mathring{\PPP}_\alpha$.

The standard bookkeeping which is used to guarantee that
$\V[G] \models \MA + 2^{\aleph_0} = \aleph_2$ 
is modified slightly here, since we need to assume
inductively that $\one\forces_{\FFF_\alpha}$ ``$(E, W)$ is good''
for all $W$ such that $(E,W)$ is good in $\V$.
This is easily seen (similarly to Theorem 49 of \cite{Je})
to be preserved at limit $\alpha$.
For the successor stage, assume that we have  $\FFF_\alpha$
and the standard bookkeeping says that we should use 
$\mathring{\QQQ}_\alpha$, which is an $\FFF_\alpha$--name
which is forced by $\one$ to be a ccc poset of size $\aleph_1$.
Roughly, we ensure that either $\MA$ holds for $\mathring{\QQQ}_\alpha$
or $\mathring{\QQQ}_\alpha$ ceases to be ccc.  More formally,
choose $\mathring\PPP_\alpha$ as follows:

Consider this from the point of view of the
$\FFF_\alpha$--extension $\V[G \cap \FFF_\alpha]$.
In this model, CH holds,
and we have a ccc poset $\QQQ_\alpha$,
and we must define another ccc poset $\PPP_\alpha$.
We know (using our inductive assumption)
that for all $W \in \V$, if $(E,W)$  good in $\V$ then it is still good.
If for all such $W$, $\one\forces_{\QQQ_\alpha}$ ``$(E, W)$ is good'',
then let $\PPP_\alpha = \QQQ_\alpha$.
If not, then fix $W \in \V$ with $(E,W)$  good in $\V$ such that
$q \forces_{\QQQ_\alpha}$ ``$(E,W)$ is not good'' for some $q \in \QQQ_\alpha$.
Still working in $\V[G \cap \FFF_\alpha]$, we apply
Lemma \ref{lemma-kill-ccc} and  let $\PPP$ be a ccc poset
of size $\aleph_1$ such that
$\QQQ_\alpha \times \PPP$ is
not ccc and such that for all $U \in \V[G \cap \FFF_\alpha]$ 
(and hence for all $U \in \V$):
If $(E,U)$ is good then $\one \forces_\PPP$ ``$(E,U)$ is good''.
Since $\QQQ_\alpha \times \PPP$ is
not ccc, we may fix $p_0 \in \PPP$ such that 
$p_0 \forces_\PPP$ ``$\QQQ_\alpha$ is not ccc''.
We cannot claim that
$\one \forces_\PPP$ ``$\QQQ_\alpha$ is not ccc'',
so let $\PPP_\alpha = p\down = \{p\in\PPP: p \le p_0\}$.
Then $\one_{\PPP_\alpha} = p
 \forces_{\PPP_\alpha}$ ``$\QQQ_\alpha$ is not ccc'',
and all good $(E,U)$ from $\V$ remain good in the $\PPP_\alpha$ extension.

Now, in $\V$, let $\mathring\PPP_\alpha$ be the name for this
$\PPP_\alpha$ as chosen above.
\end{proofof}

\begin{proofof}{Theorem \ref{thm-main}}
In the ground model $\V$, assume that
$2^{\aleph_0} = \aleph_1$ and $2^{\aleph_1} = \aleph_2$.
By CH, we may fix a (weakly) Luzin set $E \subseteq \RRR^n$ (where $n \ge 2$).
Now, apply Theorem \ref{thm-oca-fail}.
\end{proofof}

\section{Use of SOCA}
\label{sec-SOCA}

It is easily seen directly that a weakly Luzin set contradicts SOCA,
so that the Filippov space cannot exist under SOCA.
Somewhat more generally,

\begin{theorem}
\label{thm-oca}
Assume SOCA.  Let  $X$ be compact, with a 
continuous map $\pi: X \onto Y$,
where $Y$ is compact metric. 
Assume further that there is an uncountable $E \subseteq Y$
such that for $y \in E$, there are three points $x_y^i \in \pi\iv\{y\}$
for $i = 0,1,2$ and disjoint open neighborhoods $U_y^i$ of $x_y^i$
such that  $\pi(U_{y}^{i})\cap \pi(U_{y}^{j}) = \{y\}$ whenever $i \ne j$. 

Then $X$ has an uncountable discrete subset.
\end{theorem}

Note that the double arrow space satisfies these hypotheses
with ``three'' weakened to ``two'', while the Filippov space
satisfies these hypotheses with ``three'' strengthened to ``omega''.

\begin{proof}
Let $F_y^i = \cl(\pi(U_y^i))$, which is a closed
set in $Y$ containing $y$. 
Shrinking the $U_{y}^{i}$, we may assume that the three sets
$F_y^i \backslash \{y\}$ are pairwise disjoint.

We use CSM, which is a consequence of SOCA; see \cite{To}. 
Call $T \subseteq E\ $
$i$--\emph{connected} iff for all $\{y,z\} \in [T]^2$, either
$y \in F_z^i$ or $z \in F_y^i$.
Call $T$ $i$--\emph{free} iff for all $\{y,z\} \in [T]^2$, both
$y \notin F_z^i$ and $z \notin F_y^i$.
Applying CSM three times, we get an uncountable $T \subseteq E$
such that for each $i$, either
$T$ is $i$--connected or $T$ is $i$ free.   By the disjointness of the $F_y^i \backslash \{y\}$,
$T$ can be $i$--connected for at most two values of $i$.  Fixing $i$ such that
$T$ is $i$--free, we see that $\{x_y^i : y \in T\}$ is discrete.
\end{proof}


\begin{thebibliography}{99}

\bibitem{ARS}
U. Abraham, M. Rubin, and S. Shelah,
On the consistency of some partition theorems for continuous colorings,
and the structure of $\aleph\sb 1$-dense real order types,
\textit{Ann. Pure Appl. Logic} 29 (1985) 123-206.

\bibitem{AU}
P. S. Alexandroff and P. S. Urysohn,  M\'emoire sur les espaces
topologiques compacts, \textit{Verh. Akad. Wetensch. Amsterdam} 14 (1929) 1-96.

\bibitem{AS}
U. Avraham and S. Shelah,
Martin's axiom does not imply that every two $\aleph \sb{1}$-dense
sets of reals are isomorphic,
\textit{Israel J. Math.} 38 (1981) 161-176. 

\bibitem{Fil} V. V. Filippov, Perfectly normal bicompacta (Russian),
\textit{Dokl. Akad. Nauk SSSR} 189 (1969) 736-739; English translation in 
\textit{Soviet Math. Dokl.} 10 (1969) 1503-1507.

\bibitem{Gr}
G. Gruenhage,
Perfectly normal compacta, cosmic spaces, and some partition problems,
\textit{Open Problems in Topology}, North-Holland, 1990, pp. 85-95.

\bibitem{Je} T. Jech,
\textit{Lectures in Set Theory,
with Particular Emphasis on the Method of Forcing},
Springer-Verlag, 1971.

\bibitem{Ju} I. Juh\'asz,
{\it Cardinal Functions in Topology -- Ten Years Later},
Mathematical Center Tracts \#123,
Mathematisch Centrum, 1980.

\bibitem{MER} 
P. Nyikos,
Problem K.6,
\textit{Topology Proceedings} 7 (1982) 385.

\bibitem{Sze}  Z. Szentmikl\'ossy,
S-spaces and L-spaces under Martin's axiom,
\textit{Topology, Vol. II},
\textit{Colloq. Math. Soc. J\'anos Bolyai} 23,
North-Holland, 1980, pp.  1139-1145.

\bibitem{To}
S.  Todor\v cevi\'c,
\textit{Partition Problems in Topology},
Contemporary Mathematics \#84,
American Mathematical Society, 1989.

\end{thebibliography}
\end{document}